\documentclass{article}[10pt]
\usepackage{amsthm}
\usepackage{amsmath}
\usepackage{amscd}
\usepackage{amssymb}
\usepackage[all]{xy}

\newtheorem{thm}[equation]{Theorem}
\newtheorem{cor}[equation]{Corollary}

\newtheorem{lem}[equation]{Lemma}
\theoremstyle{definition}

\newtheorem{rem}[equation]{Remark}
\newtheorem{exa}[equation]{Example}

\newtheorem{que}[equation]{Question}
\numberwithin{equation}{section}

\newcommand{\ul}{\underline}

\newcommand{\surj}{\twoheadrightarrow}

\newcommand{\opn}{\operatorname}
\newcommand{\cat}[1]{\operatorname{\mathsf{#1}}}
\newcommand{\bdot}{{\textstyle \cdot}}
\newcommand{\rmitem}[1]{\item[\text{\textup{(#1)}}]}
\newcommand{\mfrak}[1]{\mathfrak{#1}}
\newcommand{\mcal}[1]{\mathcal{#1}}
\newcommand{\msf}[1]{\mathsf{#1}}

\newcommand{\mrm}[1]{\mathrm{#1}}
\newcommand{\mbb}[1]{\mathbb{#1}}

\newcommand{\tup}[1]{\textup{#1}}
\newcommand{\bsym}[1]{\boldsymbol{#1}}

\newcommand{\boplus}{\bigoplus\nolimits}

\DeclareMathSymbol{\mathbbk}{\mathord}{AMSb}{"7C}
\renewcommand{\k}{\mathbbk}

\begin{document}

\title{The Action of Adeles on the Residue Complex}
\author{Amnon Yekutieli\\
Department of Mathematics\\
Ben Gurion University\\
Be'er Sheva 84105, ISRAEL\\
amyekut@math.bgu.ac.il\\[5mm]
\textit{\small Dedicated to Professor Steven Kleiman
on the Occasion of his Sixtieth Birthday}
\setcounter{footnote}{-1}
\thanks{
Primary: 14F05; Secondary: 14F40, 13B35, 14B15.}}
\date{}
\maketitle

\begin{abstract}
Let $X$ be a scheme of finite type over a perfect field $\k$. 
In this paper we study the relation between two objects 
associated to $X$: 
the {\em Grothendieck residue complex} $\mcal{K}^{\bdot}_{X}$
and the {\em Beilinson adeles complex} 
$\ul{\mbb{A}}^{\bdot}_{\mrm{red}}(\mcal{O}_{X})$. The latter is a 
differential graded algebra (DGA).
Our first main result (Theorem \ref{thm0.1}) is that 
$\mcal{K}^{\bdot}_{X}$ is a right differential graded (DG) module over
$\ul{\mbb{A}}^{\bdot}_{\mrm{red}}(\mcal{O}_{X})$. 

We give an application to de Rham theory. Define graded sheaves
$\mcal{F}^{\bdot}_{X} := 
\mcal{H}om_{\mcal{O}_{X}}(\Omega^{\bdot}_{X / \k}, 
\mcal{K}^{\bdot}_{X})$
and
$\mcal{A}^{\bdot}_{X} := 
\ul{\mbb{A}}^{\bdot}_{\mrm{red}}(\mcal{O}_{X}) \otimes_{\mcal{O}_{X}}
\Omega^{\bdot}_{X / \k}$.
It is known that $\mcal{A}^{\bdot}_{X}$ is a DGA.
Our second main result (Theorem \ref{thm0.2}) is that 
$\mcal{F}^{\bdot}_{X}$ is a right DG $\mcal{A}^{\bdot}_{X}$-module.
When $X$ is smooth then 
$\mcal{F}^{\bdot}_{X}$ calculates de Rham homology, 
$\mcal{A}^{\bdot}_{X}$ calculates cohomology, and the action 
induces the cap product. 
We extend these constructions to singular schemes in 
characteristic $0$ using smooth formal embeddings.
\end{abstract}

\setcounter{section}{-1}
\section{Introduction}

Let $X$ be a scheme of finite type over a field $\k$. In this paper we 
study the relation between two objects associated to $X$: 
the {\em Grothendieck residue complex} $\mcal{K}^{\bdot}_{X}$
and the {\em Beilinson adeles complex} 
$\ul{\mbb{A}}^{\bdot}_{\mrm{red}}(\mcal{O}_{X})$.

Grothendieck duality theory, developed around 1960 (cf.\ \cite{RD}),
is a vast generalization of Serre duality. It is a deep and complicated 
theory, fully expressible only in the language of derived categories. 
Attempts to simplify it or find some explicit presentation of it 
attracted  
a considerable amount of research (a partial list of references is 
\cite{AK}, \cite{Kl}, \cite{Li}, \cite{HK}, \cite{KW}, \cite{LS},
\cite{Ye1}, \cite{Ne}, \cite{AJL} and \cite{Co}).

In our situation denote by $\pi : X \to \opn{Spec} \k$ the structural 
morphism. Then there is a functor
$\pi^{!} : \msf{D}^{+}_{\mrm{c}}(\cat{Mod} \k) \to 
\msf{D}^{+}_{\mrm{c}}(\cat{Mod} \mcal{O}_{X})$
between derived categories, called the twisted inverse image. The 
object $\pi^{!} \k$ is a dualizing complex on $X$. 
It has a canonical 
representative, namely its Cousin complex $\mcal{K}^{\bdot}_{X}$, 
which is called the {\em residue complex} of $X$. 
This is a bounded complex of 
quasi-coherent injective $\mcal{O}_{X}$-modules, and as a sheaf
$\mcal{K}^{\bdot}_{X} = \bigoplus_{x \in X} \mcal{K}_{X}(x)$,
where $\mcal{K}_{X}(x)$ is a constant sheaf with support 
$\overline{ \{x\} }$.
$\mcal{K}^{\bdot}_{X}$ enjoys some remarkable properties, 
that are deduced from corresponding properties of $\pi^{!}$.

Almost twenty years later Beilinson introduced his scheme theoretic
adeles (see \cite{Be}). This high dimensional generalization of the
classical adeles of a curve is actually pretty easy to define
(see Section 1). 
Given any quasi-coherent $\mcal{O}_{X}$-module
$\mcal{M}$, the complex adeles with values in $\mcal{M}$ is a 
bounded complex $\ul{\mbb{A}}^{\bdot}_{\mrm{red}}(\mcal{M})$ of 
flasque $\mcal{O}_{X}$-modules, and there is a canonical 
quasi-isomorphism 
$\mcal{M} \to \ul{\mbb{A}}^{\bdot}_{\mrm{red}}(\mcal{M})$. 
The sheaf $\ul{\mbb{A}}^{q}_{\mrm{red}}(\mcal{M})$ is a ``restricted 
product'' of local factors, each such local factor corresponding to the
geometric data of a chain $(x_{0}, \ldots, x_{q})$ of points
in $X$. Taking $\mcal{M} = \mcal{O}_{X}$ 
we obtain a DGA (differential graded algebra)
$\ul{\mbb{A}}^{\bdot}_{\mrm{red}}(\mcal{O}_{X})$.

Here is the first main result of this paper.

\begin{thm} \label{thm0.1}
Suppose $\k$ is a perfect field and $X$ is a finite type 
$\k$-scheme. Then $\mcal{K}^{\bdot}_{X}$ is a 
right DG $\ul{\mbb{A}}^{\bdot}_{\mrm{red}}(\mcal{O}_{X})$-module.
\end{thm}

The theorem is restated in more detail in Section 1 (Theorem 
\ref{thm1.1}) and proved there, using the explicit construction 
of $\mcal{K}^{\bdot}_{X}$ described in \cite{Ye3}. 
This construction is based on the theory of {\em Beilinson completion
algebras} (BCAs) developed in \cite{Ye2}.
The action of the adeles is by ``taking residues'':
multiplication by an adele supported on a chain 
$(x_{0}, \ldots, x_{q})$ is a map from $\mcal{K}_{X}(x_{0})$
to $\mcal{K}_{X}(x_{q})$.
In Question \ref{que1.1} we speculate on a generalization of Theorem 
\ref{thm0.1}.

In Section 2 we move on to de Rham theory. Let 
$(\Omega^{\bdot}_{X / \k}, \mrm{d})$
be the the algebraic de Rham complex of
$X$, which is a sheaf of commutative DGAs over $\k$. 
The graded sheaf
\[ \mcal{F}^{\bdot}_{X} := \mcal{H}om_{\mcal{O}_{X}}
(\Omega^{\bdot}_{X / \k}, \mcal{K}^{\bdot}_{X}) \]
is a graded $\Omega^{\bdot}_{X / \k}$-module. 
According to \cite{Ye3} $\mcal{F}^{\bdot}_{X}$ has a 
coboundary operator $\mrm{D}$ that's a differential operator 
of order $\leq 1$ over $\mcal{O}_{X}$, and
$(\mcal{F}^{\bdot}_{X}, \mrm{D})$ is a DG 
$\Omega^{\bdot}_{X / \k}$-module.

On the other hand we have the graded algebra
\[ \mcal{A}^{\bdot}_{X} := 
\boplus_{p, q} \
\ul{\mbb{A}}^{q}_{\mrm{red}}(\Omega^{p}_{X / \k})
\cong
\ul{\mbb{A}}^{\bdot}_{\mrm{red}}
(\mcal{O}_{X}) \otimes_{\mcal{O}_{X}} \Omega^{\bdot}_{X / \k} . \]
By \cite{HY1} there is a differential $\mrm{D}$ that makes
$(\mcal{A}^{\bdot}_{X}, \mrm{D})$ into a DGA, and moreover
$\Omega^{\bdot}_{X / \k} \to \mcal{A}^{\bdot}_{X}$
is a DGA quasi-isomorphism. Theorem \ref{thm0.1} implies that 
$\mcal{F}^{\bdot}_{X}$ is a graded right 
$\mcal{A}^{\bdot}_{X}$-module.
Our second main result is:

\begin{thm} \label{thm0.2}
Suppose $\k$ is a perfect field and $X$ is a finite type 
$\k$-scheme. Then $\mcal{F}^{\bdot}_{X}$
is a DG right $\mcal{A}^{\bdot}_{X}$-module.
\end{thm}

Theorem \ref{thm0.2} is used in \cite{HY2} to state the adelic 
Gauss-Bonnet formula, which is proved there. See also Remark 
\ref{rem2.3} below.

Consider a finite type $\k$-scheme $Y$ and a closed subscheme 
$X \subset Y$. Define $\mfrak{X} := Y_{/ X}$, the formal 
completion of $Y$ along $X$ (cf.\ \cite{EGA} I).
Let $X_i$ be the $i$th infinitesimal neighborhood of 
$X = X_0$ in $Y$. Define
$\mcal{A}^{\bdot}_{\mfrak{X}} := \lim_{\leftarrow i}
\mcal{A}^{\bdot}_{X_i}$
and
$\mcal{F}^{\bdot}_{\mfrak{X}} := \lim_{i \to}
\mcal{F}^{\bdot}_{X_i}$.

\begin{cor} \label{cor0.1}
$\mcal{F}^{\bdot}_{\mfrak{X}}$ is a right DG 
$\mcal{A}^{\bdot}_{\mfrak{X}}$-module. 
If $\opn{char} \k = 0$ and $Y$ is smooth over $\k$
then de Rham cohomology is
$\mrm{H}^p_{\mrm{DR}}(X) = 
\mrm{H}^p \Gamma(X, \mcal{A}^{\bdot}_{\mfrak{X}})$,
de Rham homology is 
$\mrm{H}_p^{\mrm{DR}}(X) = 
\mrm{H}^{-p} \Gamma(X, \mcal{F}^{\bdot}_{\mfrak{X}})$,
and the right action of $\mcal{A}^{\bdot}_{\mfrak{X}}$
on $\mcal{F}^{\bdot}_{\mfrak{X}}$ induces the cap product.
\end{cor}

The embedding $X \subset \mfrak{X}$ is an instance of a smooth 
formal embedding in the sense of \cite{Ye4}.

Theorem \ref{thm0.2} and Corollary \ref{cor0.1} are proved at the
end of Section 2.

\medskip \noindent
{\em Acknowledgments.}\
The paper is dedicated to Steven Kleiman in gratitude for his 
longtime interest and support of my work on residues. 
I also wish to thank Reinhold H\"{u}bl, Joseph Lipman and 
Pramathanath Sastry for helpful discussions. Finally I wish to 
express thanks to the referee for reading the paper carefully and 
suggesting improvements.


\section{The Action}

Let us begin with a review of {\em Beilinson adeles} on a noetherian
scheme $X$. A chain of length $q$ of points in $X$ is a sequence 
$\xi = (x_{0}, \ldots, x_{q})$ of points with
$x_{i+1} \in \overline{ \{ x_{i} \} }$. Denote by $S(X)_{q}$ the set
of length $q$ chains, so $\{ S(X)_{q} \}_{q \geq 0}$ is a simplicial 
set. For a subset $T \subset S(X)_{q}$ and a point $x \in X$ let
\[ \hat{x} T := \{ (x_{1}, \ldots, x_{q}) \mid (x, x_{1}, \ldots, x_{q})
\in T \} . \] 
Denote by $\mfrak{m}_x \subset \mcal{O}_{X, x}$ the maximal ideal.
If $\mcal{M}$ is a coherent $\mcal{O}_{X}$-module then for any 
$n \geq 1$ the $\mcal{O}_{X, x}$-module 
$\mcal{M}_{x} / \mfrak{m}_{x}^{n} \mcal{M}_{x}$ 
can be thought of as a quasi-coherent sheaf, constant on the 
closed set $\overline{ \{ x \} }$.
According to \cite{Be} there is a unique collection of functors
$\mbb{A}(T, -) : \cat{QCoh} \mcal{O}_{X} \to \cat{Ab}$,
indexed by subsets $T \subset S(X)_{q}$, 
each of which commuting with direct limits, and satisfying 
\begin{equation} \label{eqn1.10}
\mbb{A}(T, \mcal{M}) = \begin{cases}
\prod_{(x) \in T} \lim_{\leftarrow n} \mcal{M}_{x} / \mfrak{m}_{x}^{n}
\mcal{M}_{x} & \text{if } q = 0 \\[1em]
\prod_{x \in X} \lim_{\leftarrow n} 
\mbb{A}(\hat{x} T, \mcal{M}_{x} / \mfrak{m}_{x}^{n} \mcal{M}_{x})
& \text{if } q > 0 
\end{cases} 
\end{equation}
for $\mcal{M}$ coherent.
Furthermore each $\mbb{A}(T, -)$ is exact.

For a single chain $\xi$ one also writes
$\mcal{M}_{\xi} := \mbb{A}(\{ \xi \}, \mcal{M})$, 
and this is the {\em Beilinson completion} of $\mcal{M}$ along $\xi$.
Then
\begin{equation} \label{eqn1.6}
\mbb{A}(T, \mcal{M}) \subset \prod_{\xi \in T} \mcal{M}_{\xi} .
\end{equation}
In view of this we shall say that $\mbb{A}(T, \mcal{M})$ is the 
group of adeles combinatorially supported on $T$ and with values in
$\mcal{M}$.

Observe that for $q = 0$ and $\mcal{M}$ coherent we have
$\mcal{M}_{(x)} = \widehat{\mcal{M}}_{x}$, the $\mfrak{m}_{x}$-adic
completion, and (\ref{eqn1.6}) is an equality.

Define a presheaf 
$\ul{\mbb{A}}(T, \mcal{M})$ by 
\begin{equation} \label{eqn1.4}
\Gamma(U, \ul{\mbb{A}}(T, \mcal{M})) := 
\mbb{A}(T \cap S(U)_{q}, \mcal{M}) 
\end{equation}
for $U \subset X$ open. Then $\ul{\mbb{A}}(T, \mcal{M})$ is in 
fact a flasque sheaf. Also $\ul{\mbb{A}}(T, \mcal{O}_{X})$ is a flat 
$\mcal{O}_{X}$-algebra, and
$\ul{\mbb{A}}(T, \mcal{M}) \cong \ul{\mbb{A}}(T, \mcal{O}_{X})
\otimes_{\mcal{O}_{X}} \mcal{M}$. 
Given a local section 
$m \in \ul{\mbb{A}}(T, \mcal{M})$ we shall often use the inclusion
(\ref{eqn1.6}) to write
$m = (m_{\xi})$ where $\xi$ runs over $T$ and 
$m_{\xi} \in \mcal{M}_{\xi}$. 

Let $S(X)_{q}^{\mrm{red}}$ be the set of reduced chains
(i.e.\ without repeated points), and define
\[ \ul{\mbb{A}}^{q}_{\mrm{red}}(\mcal{M}) :=
\ul{\mbb{A}}(S(X)_{q}^{\mrm{red}}, \mcal{M}) . \]
For $0 \leq i \leq q$ the $i$th face map $\partial_i$, 
which omits the point 
$x_i$ from a chain $(x_0, \ldots, x_q)$, induces a homomorphism 
\[ \partial^i: \ul{\mbb{A}}^{q - 1}_{\mrm{red}}(\mcal{M}) \to
\ul{\mbb{A}}^{q}_{\mrm{red}}(\mcal{M}) . \]
Then $\ul{\mbb{A}}^{\bdot}_{\mrm{red}}(\mcal{M})$ 
is a complex with coboundary operator
$\partial := \sum (-1)^{i} \partial^{i}$,
and 
$\mcal{M} \to \ul{\mbb{A}}^{\bdot}_{\mrm{red}}(\mcal{M})$
is a quasi-isomorphism. Since 
$\ul{\mbb{A}}^{\bdot}_{\mrm{red}}(\mcal{M})$
is a complex of flasque sheaves we get
\[ \mrm{H}^{q} \Gamma(X, \ul{\mbb{A}}^{\bdot}_{\mrm{red}}(\mcal{M}))
= \mrm{H}^{q}(X, \mcal{M}) . \]

The complex 
$\ul{\mbb{A}}^{\bdot}_{\mrm{red}}(\mcal{O}_{X})$
is a DGA, with the Alexander-Whitney product. For local sections
$a \in \ul{\mbb{A}}^{q}_{\mrm{red}}(\mcal{O}_{X})$
and
$b \in \ul{\mbb{A}}^{q'}_{\mrm{red}}(\mcal{O}_{X})$
the product is
\[ a \cdot b := \partial^{-}(a) \cdot \partial^{+}(b) \in
\ul{\mbb{A}}^{q+q'}_{\mrm{red}}(\mcal{O}_{X}) , \]
where $\partial^{-}$ and $\partial^{+}$ correspond respectively to the 
initial and final segments of $(0, \ldots, q, \ldots, q+q')$.
This algebra is not (graded) commutative. 
For proofs and more details see \cite{Hr}, \cite{Ye1} Chapter 3
and \cite{HY1} Section 1.

\begin{exa} \label{exa1.2}
Suppose $X$ is a nonsingular curve. The relation to the classical 
ring of adeles $\mbb{A}(X)$ of Chevalley and Weil associated to $X$
is
\[ \mbb{A}(X) = 
\Gamma(X, \ul{\mbb{A}}^{1}_{\mathrm{red}}(\mcal{O}_{X})) . \]
\end{exa}

 From now we assume $X$ is a finite type scheme over a perfect field
$\k$.

Next let us recall the construction of the residue complex
$\mcal{K}^{\bdot}_{X}$ in \cite{Ye3}.
It starts with the theory of Beilinson completion algebras (BCAs)
developed in \cite{Ye2}. A BCA $A$ is a semilocal $\k$-algebra with a 
topology and with valuations on its residue fields. Each local factor
of $A$ is a quotient of 
$K((s_{1}, \ldots, s_{n}))[[t_{1}, \ldots, t_{m}]]$
where $K$ is some finitely generated extension field of $\k$,
and $K((s_{1}, \ldots, s_{n})) = K((s_{n})) \cdots ((s_{1}))$
is the field of iterated Laurent series.
One considers two kinds of homomorphisms between BCAs:
morphisms $f : A \to B$ and intensifications $u : A \to \widehat{A}$.

Each BCA $A$ has a dual module $\mcal{K}(A)$, which is functorial 
w.r.t.\ these homomorphisms; namely there are maps
$\opn{Tr}_{f} : \mcal{K}(B) \to \mcal{K}(A)$
and
$\mrm{q}_{u} : \mcal{K}(A) \to \mcal{K}(\widehat{A})$.
Algebraically $\mcal{K}(A)$ is an injective hull of $A / \mfrak{r}$, 
where $\mfrak{r}$ is the Jacobson radical. 

\begin{exa}
Let $A := \k [[s_1, s_2]]$, $B := \k (s_1, s_2)$ and 
$\widehat{B} := \k ((s_1, s_2))$. 
The inclusions
$f: A \to \widehat{B}$ and $u: B \to \widehat{B}$
are respectively a morphism and an intensification.
The dual modules are
\[ \mcal{K}(A) = \opn{Hom}^{\mrm{cont}}_{\k}(A, \k) , \]
\[ \mcal{K}(B) = \Omega^2_{B / \k} = 
B \cdot \mrm{d} s_1 \wedge \mrm{d} s_2  \]
and
\[ \mcal{K}(\widehat{B}) = 
\Omega^{2, \mrm{sep}}_{\widehat{B} / \k} = 
\widehat{B} \cdot \mrm{d} s_1 \wedge \mrm{d} s_2 . \]
The homomorphism 
$\mrm{q}_u: \mcal{K}(B) \to \mcal{K}(\widehat{B})$ 
is the inclusion. For a form 
$\beta \in \Omega^{2, \mrm{sep}}_{\widehat{B} / \k}$
the functional $\opn{Tr}_f (\beta) \in \mcal{K}(A)$ is described 
as follows. Given an element $a \in A$ write
\[ a \beta = \sum_{i, j} \lambda_{i,j} s_1^i s_2^j \cdot
\mrm{d} s_1 \wedge \mrm{d} s_2  \in 
\Omega^{2, \mrm{sep}}_{\widehat{B} / \k} \]
with $\lambda_{i,j} \in \k$. Then 
\[ \opn{Tr}_f (\beta)(a) = \lambda_{-1, -1} , \]
namely the residue of $a \beta$.
\end{exa}

Suppose $\xi = (x, \ldots, y)$ is a saturated chain of points in
$X$ (i.e.\ each point is an immediate specialization of the 
previous one). Then the Beilinson completion 
$\mcal{O}_{X, \xi}$ is a BCA. The natural algebra homomorphisms
$\partial^{-} : \mcal{O}_{X, (x)} \to \mcal{O}_{X, \xi}$ and
$\partial^{+} : \mcal{O}_{X, (y)} \to \mcal{O}_{X, \xi}$
are an intensification and a morphism, respectively. So there are 
homomorphisms on dual modules
\[ \mrm{q}_{\partial^{-}} : \mcal{K}(\mcal{O}_{X, (x)}) 
\to \mcal{K}(\mcal{O}_{X, \xi}) \]
and
\[ \opn{Tr}_{\partial^{+}} : \mcal{K}(\mcal{O}_{X, \xi}) \to 
\mcal{K}(\mcal{O}_{X, (y)}) . \]
The composition 
\[ \opn{Tr}_{\partial^{+}} \circ\ \mrm{q}_{\partial^{-}}:
\mcal{K}(\mcal{O}_{X, (x)}) \to \mcal{K}(\mcal{O}_{X, (y)}) \]
is denoted by $\delta_{\xi}$. We regard 
$\mcal{K}_{X}(x) := \mcal{K}(\mcal{O}_{X, (x)})$ 
as a quasi-coherent 
$\mcal{O}_{X}$-module, constant on the closed set $\overline{\{ x \}}$.
Define
\begin{equation} \label{eqn1.1}
 \mcal{K}_{X}^{-q} := \bigoplus_{\opn{dim} \overline{\{x\}} = q}
\mcal{K}_{X}(x) 
\end{equation}
and 
\begin{equation} \label{eqn1.2} 
\delta := (-1)^{q+1} \sum_{(x,y)} \delta_{(x,y)} :
\mcal{K}_{X}^{-q} \to \mcal{K}_{X}^{-q + 1} . 
\end{equation}
Then the pair $(\mcal{K}^{\bdot}_{X}, \delta)$ is the residue complex of
$X$. That is to say, in the notation of the introduction,
there is a canonical isomorphism
$\mcal{K}^{\bdot}_{X} \cong \pi^{!} \k$ in the derived
category $\msf{D}(\cat{Mod} \mcal{O}_{X})$ (see \cite{Ye3} Corollary 2.5).

Let $x$ be a point of dimension $q$ in $X$, and consider a local section
$\phi_{x} \in \mcal{K}_{X}(x) \subset \mcal{K}_{X}^{-q}$. 
Let $\xi = (x_{0}, \ldots, x_{q'})$ be any chain of length $q'$ in $X$, 
and let
$a_{\xi} \in \mcal{O}_{X, \xi}$. 
Define an element
$\phi_{x} \cdot a_{\xi} \in \mcal{K}_{X}^{-q + q'}$ 
as follows. 
If $x = x_{0}$ and $\xi$ is saturated then 
there are homomorphisms
$\mrm{q}_{\partial^{-}} : \mcal{K}_{X}(x) \to 
\mcal{K}(\mcal{O}_{X, \xi})$
and
$\opn{Tr}_{\partial^{+}} : \mcal{K}(\mcal{O}_{X, \xi}) \to
\mcal{K}_{X}(x_{q'})$.
Since $\mcal{K}(\mcal{O}_{X, \xi})$ is an 
$\mcal{O}_{X,\xi}$-module the product
$a_{\xi} \cdot \mrm{q}_{\partial^{-}} (\phi_{x}) \in
\mcal{K}(\mcal{O}_{X, \xi})$ 
exists, and we set
\begin{equation} \label{eqn1.9}
\phi_{x} \cdot a_{\xi} := 
\opn{Tr}_{\partial^{+}}(a_{\xi} \cdot \mrm{q}_{\partial^{-}} (\phi_{x}))
\in \mcal{K}_{X}(x_{q'}) .
\end{equation}
Otherwise we set
$\phi_{x} \cdot a_{\xi} := 0$.

\begin{rem}
In order to apply the Koszul sign rule we consider the ``local'' 
objects $\mcal{K}_{X}(x)$ and $\mcal{O}_{X, \xi}$ as ungraded; 
whereas the ``global'' objects $\mcal{K}^{\bdot}_{X}$ and
$\ul{\mbb{A}}^{\bdot}_{\mrm{red}}(\mcal{O}_{X})$
are considered as graded. Ungraded local elements shall be 
decorated with suitable subscripts, such as 
$a_{\xi} \in \mcal{O}_{X, \xi}$ or
$\phi_x \in \mcal{K}_{X}(x)$.
\end{rem}

Suppose $\xi = (x_0, \ldots, x_p)$ and 
$\eta = (y_0, \ldots, y_q)$ are chains such that 
$y_0 \in \overline{\{ x_p \}}$. Then we denote by
$\xi \vee \eta$ the concatenated chain 
$(x_0, \ldots, x_p, y_0, \ldots, y_q)$. 

\begin{lem} \label{lem1.1}
Let $\xi = (x, \ldots, y)$ and $\eta = (y, \ldots, z)$
be saturated chains in $X$ and let 
$\phi_x \in \mcal{K}_{X}(x)$, 
$a_\xi \in \mcal{O}_{X, \xi}$ and
$b_\eta \in \mcal{O}_{X, \eta}$ be some elements. Writing 
$\xi \vee \partial_{0}(\eta) := (x, \ldots, y, \ldots, z)$, one has
$a_\xi \cdot b_\eta \in \mcal{O}_{X, \xi \vee \partial_{0} (\eta)}$. 
Then
\[ (\phi_x \cdot a_\xi) \cdot b_\eta = 
\phi_x \cdot (a_\xi \cdot b_\eta) \in 
\mcal{K}_{X}(z). \]
\end{lem}

\begin{proof}
Consider the commutative diagram with intensifications 
$u, v, w, u'$ and morphisms $f, f', g, h$:
\[ \UseTips \xymatrix{ 
\mcal{O}_{X, (x)} \ar[r]^{w}  \ar[dr]_{v} 
& \mcal{O}_{X, \xi} \ar[d]^{u'}
& \mcal{O}_{X, (y)} \ar[l]_{f} \ar[d]^{u} \\
& \mcal{O}_{X, \xi \vee \partial_{0} (\eta)}
& \mcal{O}_{X, \eta} \ar[l]_{f'} \\
&  & \mcal{O}_{X, (z)} \ar[ul]^{h} \ar[u]_{g} 
} \]
According to \cite{Ye3} Lemma 1.12 there is an isomorphism
\[ \mcal{O}_{X, \xi \vee \partial_{0} (\eta)} \cong
\mcal{O}_{X, \xi} \otimes^{(\wedge)}_{\mcal{O}_{X, (y)}}
\mcal{O}_{X, \eta} \]
(intensification base change), so by
\cite{Ye2} Thm.\ 7.4 (ii) one has
\[ \opn{Tr}_{f'} \circ\ \mrm{q}_{u'}  
= \mrm{q}_{u} \circ\ \opn{Tr}_{f} . \] 

Now the BCA $\mcal{O}_{X, \xi \vee \partial_{0} (\eta)}$ is 
commutative, $\opn{Tr}_{f'}$ is $\mcal{O}_{X, \eta}$-linear 
and $\mrm{q}_{u'}$ is $\mcal{O}_{X, \xi}$-linear. Hence
\[ \begin{aligned}
\phi_x \cdot (a_\xi \cdot b_\eta) & = 
\opn{Tr}_{h} \bigl( (a_\xi \cdot b_\eta) \cdot 
\mrm{q}_{v}(\phi_x) \bigr) \\
& = \opn{Tr}_{g} \bigl( \opn{Tr}_{f'} \bigl(
f'(b_\eta) \cdot u'(a_\xi) \cdot
\mrm{q}_{u'}(\mrm{q}_{w}(\phi_x)) \bigr) \bigr) \\
& = \opn{Tr}_{g} \bigl( b_\eta \cdot \opn{Tr}_{f'}
\bigl( \mrm{q}_{u'}(a_\xi \cdot \mrm{q}_{w}(\phi_x)) \bigr) \bigr) \\
& = \opn{Tr}_{g} \bigl( b_\eta \cdot \mrm{q}_{u}
\bigl( \opn{Tr}_{f}(a_\xi \cdot \mrm{q}_{w}(\phi_x)) \bigr) \bigr) \\
& = (\phi_x \cdot a_\xi) \cdot b_\eta .
\end{aligned} \]
\end{proof}

Suppose $x \in X$ is a point and $b \in \mcal{O}_{X, (x)}$ 
is some element. Using the inclusion
$\mcal{O}_{X, (x)} \subset \ul{\mbb{A}}_{\mrm{red}}^0(\mcal{O}_{X})$
we consider $b$ as an adele.

\begin{lem}[Approximation] \label{lem1.4}
Let $U \subset X$ be an affine open set, $x \in U$ a point 
and $q \geq 1$. Define
\[ T := \{ \xi \in S(U)^{\mrm{red}}_q \mid \xi = (x, \ldots) \} . 
\]
Let $a \in \mbb{A}(T, \mcal{O}_{X})$ be some adele. 
\begin{enumerate}
\rmitem{1} Given an integer $n \geq 1$ there exists adeles 
$b \in \mcal{O}_{X, (x)}$ and
$c = (c_{\eta}) \in \mbb{A}(\partial_0(T), \mcal{O}_{X})$
such that
\[ a - \partial^1(b) \cdot c \in \mfrak{m}_x^n \cdot
\mbb{A}(T, \mcal{O}_{X}) . \]
\rmitem{2} Given an element $\phi_x \in \mcal{K}_X(x)$ there is an 
integer $n \geq 1$ such that 
$\mfrak{m}_x^n \cdot \phi_x = 0$. For such $n$ the adeles $b$ and 
$c$ from part \tup{(1)} satisfy 
\[ \phi_x \cdot a_{\xi} = \phi_x \cdot (b_{(x, y)} \cdot c_{\eta})
\in \mcal{K}_X(z) \]
for all $\xi \in T$, where we write
$\eta = (y, \ldots, z) := \partial_0(\xi)$, and
$b_{(x, y)}$ is the $(x, y)$ component of $\partial^1(b)$. 
\end{enumerate}
\end{lem}

\begin{proof}
(1) There is a ring homomorphism
$\mcal{O}_{X, (x)} \to \mbb{A}(T, \mcal{O}_{X})$
and via this homomorphism we obtain an ideal
$\mfrak{m}_x^n \cdot \mbb{A}(T, \mcal{O}_{X})$. 
Because the functor $\mbb{A}(T, -)$ is exact we have
\[ \frac{\mbb{A}(T, \mcal{O}_{X})}
{\mfrak{m}_x^n \cdot \mbb{A}(T, \mcal{O}_{X})}
\cong \mbb{A}(T, \mcal{O}_{X, x} / \mfrak{m}_x^n) . \]
By the definition of adeles with values in a quasi-coherent 
sheaf there is an isomorphism
\[ \mbb{A}(T, \mcal{O}_{X, x} / \mfrak{m}_x^n)
\cong \lim_{\to} \mbb{A}(T, \mcal{L}) \]
where the limit runs over the coherent $\mcal{O}_{U}$-submodules 
$\mcal{L} \subset \mcal{O}_{X, x} / \mfrak{m}_x^n$.
We note that $\partial_0(T) = \hat{x} T$, so for $\mcal{L}$ large 
enough (i.e.\ $\mcal{L}_x = \mcal{O}_{X, x} / \mfrak{m}_x^n$)
\[ \begin{aligned}
\mbb{A}(T, \mcal{L}) & = 
\opn{lim}_{\leftarrow i} 
\mbb{A}(\partial_0(T), \mcal{L}_x / \mfrak{m}_x^i \mcal{L}_x) \\
& \cong 
\mbb{A}(\partial_0(T), \mcal{O}_{X, x} / \mfrak{m}_x^n) .
\end{aligned} \]
The conclusion is that 
\[ \mbb{A}(T, \mcal{O}_{X, x} / \mfrak{m}_x^n) \cong
\mbb{A}(\partial_0(T), \mcal{O}_{X, x} / \mfrak{m}_x^n) . \]
Again we go to coherent subsheaves.
Write $C := \Gamma(U, \mcal{O}_{X})$ and
$\mcal{O}_Z := \opn{Im}(\mcal{O}_{X}
\to \mcal{O}_{X, x} / \mfrak{m}_x^n)$. 
Then
\[ \begin{aligned}
\mbb{A}(\partial_0(T), \mcal{O}_{X, x} / \mfrak{m}_x^n)
& = \opn{lim}_{\to} \mbb{A}(\partial_0(T), \mcal{L}) \\
& \cong (\mcal{O}_{X, x} / \mfrak{m}_x^n) \otimes_C
\mbb{A}(\partial_0(T), \mcal{O}_Z) \\
& \cong (\mcal{O}_{X, x} / \mfrak{m}_x^n) \otimes_C
\mbb{A}(\partial_0(T), \mcal{O}_X) . 
\end{aligned} \]
So there is a ring surjection
\begin{equation} \label{eqn1.5}
\mcal{O}_{X, x}  \otimes_C
\mbb{A}(\partial_0(T), \mcal{O}_X) \surj
\mbb{A}(T, \mcal{O}_{X, x} / \mfrak{m}_x^n) .
\end{equation}

Consider the image $\bar{a}$ of the adele $a$ in 
$\mbb{A}(T, \mcal{O}_{X, x} / \mfrak{m}_x^n) $. 
Using the surjection (\ref{eqn1.5}) we can write
$\bar{a} = \sum_{i = 1}^r b_i \otimes c_i$
with $b_i \in \mcal{O}_{X, x}$ and 
$c_i \in \mbb{A}(\partial_0(T), \mcal{O}_X)$.
By bringing the $b_i$ to a common denominator we can assume $r = 
1$.

\medskip \noindent
(2) Because $\mcal{K}(\mcal{O}_{X, (x)})$ is an 
$\mfrak{m}_x$-torsion  module the element $\phi_x$ is 
annihilated by some power 
$\mfrak{m}_x^n$. Pick adeles $b$ and $c$ as in part (1). Then for 
any chain 
$\xi = (x) \vee \eta = (x, y, \ldots) \in T$ we obtain
\[ a_{\xi} - b_{(x, y)} \cdot c_{\eta} \in 
\mfrak{m}_x^n \cdot \mcal{O}_{X, \xi} . \]
Therefore
\[ \phi_x \cdot a_{\xi} = \phi_x \cdot (b_{(x, y)} \cdot c_{\eta}) 
. \]
\end{proof}

\begin{lem} \label{lem1.2}
Let $U \subset X$ be an open set, $x \in U$ a point, 
$\phi_x \in \mcal{K}_X(x)$ an element and 
$a = (a_{\xi}) \in \mbb{A}(S(U)^{\mrm{red}}_{q},\mcal{O}_{X})$
an adele. Then for all but finitely many chains 
$\xi \in S(U)_{q}^{\mrm{red}}$ 
one has $\phi_x \cdot a_{\xi} = 0$.
\end{lem}

\begin{proof} 
The proof is by induction on $q$. For $q=0$ there is nothing 
to prove, so take $q \geq 1$. 
Because $U$ is covered by finitely many affine open sets, and we 
are only interested in establishing finiteness, we might as well 
assume $U$ itself is affine.
Moreover by the definition of the product we can neglect those 
chains in the combinatorial support of $a$ that do not begin with 
$x$. Thus we can assume 
$a \in \mbb{A}(T, \mcal{O}_X)$,
where $T$ is the set defined in Lemma \ref{lem1.4}.
According to this lemma we can find adeles 
$b \in \mcal{O}_{X, (x)}$ and
$c = (c_{\eta}) \in \mbb{A}(\partial_0(T), \mcal{O}_{X})$
such that
\[ \phi_x \cdot a_{\xi} = \phi_x \cdot (b_{(x, y)} \cdot c_{\eta})
\in \mcal{K}_X(z) \]
for all $\xi = (x) \vee \eta = (x, y, \ldots, z) \in T$.

The only way to get a nonzero 
product $\phi_x \cdot a_\xi$ is when $\xi$
is a saturated chain. Consider such a chain $\xi$. According to
Lemma \ref{lem1.1} we have
\[ \phi_x \cdot (b_{(x, y)} \cdot c_\eta) 
= (\phi_x \cdot b_{(x, y)}) \cdot c_\eta  . \]
For any point $y$ occurring one has
\[ \phi_x \cdot b_{(x, y)} = 
\delta_{(x, y)}(b \phi_x) \in \mcal{K}_X(y) . \]
It follows that the product $\psi_y := \phi_x \cdot b_{(x, y)}$
vanishes for all but finitely many 
points $y$. Fixing $y$, the induction hypothesis applied to the 
element $\psi_y \in \mcal{K}_X(y)$
says that $\psi_y \cdot c_\eta = 0$
for all but finitely many chains 
$\eta \in S(U)_{q - 1}^{\mrm{red}}$.
\end{proof}

\begin{thm} \label{thm1.1}
$\mcal{K}_{X}^{\bdot}$ is a right DG 
$\ul{\mbb{A}}^{\bdot}_{\mrm{red}}(\mcal{O}_{X})$-module, 
with product
\[ \phi \cdot a = 
\sum_{x, \xi} \phi_{x} \cdot a_{\xi}  \]
for local sections 
$\phi = \sum_{x} \phi_{x} \in \mcal{K}_{X}^{\bdot}$ 
and
$a = (a_{\xi}) \in \ul{\mbb{A}}_{\mrm{red}}^{\bdot}(\mcal{O}_{X})$.
\end{thm}

\begin{proof}
According to Lemmas \ref{lem1.2} and \ref{lem1.1} 
this is a well defined associative product.
It remains to verify that 
\begin{equation} \label{eqn1.3}
 \delta(\phi \cdot a) = \delta (\phi) \cdot a + 
(-1)^{q} \phi \cdot \partial (a) 
\end{equation}
for 
$\phi \in \mcal{K}_{X}^{-q}$ 
and
$a \in \ul{\mbb{A}}_{\mrm{red}}^{q'}(\mcal{O}_{X})$.
We may assume $\phi = \phi_{x}$ and $a = a_{\xi}$
for some point $x$ of dimension $q$ and a chain
$\xi = (x_{0}, \ldots, x_{q'})$.

Now there are only 3 ways to get any nonzero term in equation 
(\ref{eqn1.3}):
(i) $\xi$ is saturated and $x_{0} = x$; 
(ii) $\xi$ is saturated and $(x, x_{0})$ is saturated;
or 
(iii) $x_{0} = x$, and for some index
$0 \leq i < q'$ and some point $y \in X$, the chain
\[ \eta := (x_{0}, \ldots, x_{i}, y, x_{i+1}, \ldots, x_{q'}) \]
is saturated.

In case (i),
$\delta (\phi) \cdot a = 0$ and
\[ \begin{aligned}
\phi \cdot \partial (a) 
& = \phi_x \cdot (-1)^{q'+1} \partial^{q'+1} (a_\xi) \\
& = (-1)^{q'+1} \sum_{y} \delta_{(x_{q'}, y)} (\phi_x \cdot a_\xi) 
\\
& = (-1)^{q} \delta (\phi \cdot a) . 
\end{aligned} \]

In case (ii), 
$\phi \cdot a = 0$ and
\[ \begin{aligned}
\phi \cdot \partial (a) & = 
\phi_x \cdot \partial^{0} (a)_{(x) \vee \xi} \\
& = \delta_{(x, x_{0})} (\phi_x) \cdot a_\xi \\
& = (-1)^{q + 1} \delta (\phi) \cdot a .
\end{aligned} \]
In this equation
$(x) \vee \xi$ is the concatenated chain $(x, x_0, \ldots, x_{q'})$.

Finally in case (iii), 
$\delta(\phi) \cdot a = 0$, $\delta(\phi \cdot a) = 0$, 
and it remains to show that also $\phi \cdot \partial(a) = 0$.
We note that 
\[ \phi \cdot \partial(a) = 
\sum_y \phi_x \cdot (-1)^{i + 1} \partial^{i + 1}(a)_{\eta} \]
where $y$ runs over the points such that 
$(x_i, y, x_{i + 1})$ is a saturated chain, and $\eta$ is as 
above (and depends on $y$).

For any index $0 \leq j \leq q'$ let us write
$\xi_j := (x_j, \ldots, x_{q'})$.
We shall use an approximation trick of Lemma \ref{lem1.4}
 to define recursively elements
$\phi_{x_{j}} \in \mcal{K}_X(x_{j})$ and
$b_{(x_{j})} \in \mcal{O}_{X, (x_{j})}$ 
for $0 \leq j \leq i$, and 
$a_{\xi_{j}} \in \mcal{O}_{X, \xi_{j}}$
for $0 \leq j \leq i + 1$. 
Pick some affine open set $U \subset X$ containing $x_{q'}$
($U$ does not play an essential role, since our computation is local 
anyhow; but it appears in Lemma \ref{lem1.4}). 
For $j = 0$ we note that $x_0 = x$ and $\xi_0 = \xi$, 
and the elements $\phi_{x_0}$ and $a_{\xi_0}$ are
already defined. Now suppose $j \leq i$ and 
$\phi_{x_{j}}$ and $a_{\xi_{j}}$ have been defined. 
By Lemma \ref{lem1.4} we can find elements 
$b_{(x_{j})} \in \mcal{O}_{X, (x_{j})}$ and 
$a_{\xi_{j + 1}} \in \mcal{O}_{X, \xi_{j + 1}}$ 
such that 
\[ \phi_{x_{j}} \cdot a_{\xi_{j}} = 
\phi_{x_{j}} \cdot (\partial^1(b_{(x_{j})})_{(x_j, x_{j + 1})} 
\cdot a_{\xi_{j + 1}}) \in \mcal{K}(x_{q'}) . \]

If $j < i$ we also define
\[ \phi_{x_{j + 1}} := 
\phi_{x_{j}} \cdot \partial^1(b_{(x_{j})})_{(x_j, x_{j + 1})} 
\in \mcal{K}_X(x_{j + 1}) . \]
We thus have for any of the points $y$ under consideration:
\[ \phi_{x_{j}} \cdot \partial^{(i + 1) - j}(a_{\xi_{j}})_{\eta_j}
= \phi_{x_{j + 1}} \cdot \partial^{(i + 1) - (j + 1)}
(a_{\xi_{j + 1}})_{\eta_{j + 1}} , \]
where 
\[ \eta_{j} := (x_{j}, \ldots, x_i, y, x_{i + 1}, \ldots, x_{q'}) 
. \]
Putting it all together we obtain
\[ \begin{aligned}
\phi_{x_{0}} \cdot \partial^{i + 1}(a_{\xi_{0}})_{\eta_0}
& = \phi_{x_{i}} \cdot \partial^{1}(a_{\xi_{i}})_{\eta_i} \\
& =  (\phi_{x_{i}} \cdot b_{(x_i)})
\cdot (\partial^{0} \circ \partial^{0}) (a_{\xi_{i + 1}})_{\eta_i} \\
& =  \delta_{(x_i, y, x_{i + 1})}
(\phi_{x_{i}} \cdot b_{(x_i)}) \cdot a_{\xi_{i + 1}} , 
\end{aligned} \]
and hence
\[ \phi \cdot \partial(a) 
= (-1)^{i + 1} \sum_y \delta_{(x_i, y, x_{i + 1})}
(\phi_{x_{i}} \cdot b_{(x_i)}) \cdot a_{\xi_{i + 1}} . \]
But according to \cite{Ye3} Lemma 2.15(3), which is a variant
of the Parshin-Lomadze Residue Theorem, we have
\[ \sum_y \delta_{(x_i, y, x_{i + 1})}
(\phi_{x_{i}} \cdot b_{(x_i)}) = 0 .\]
\end{proof}

\begin{exa} \label{exa1.1}
The global section 
$a := \partial^1(1) \in 
\Gamma(X, \ul{\mbb{A}}^{1}_{\mrm{red}}(\mcal{O}_{X}))$
acts on $\mcal{K}_{X}^{q}$ like $(-1)^{q + 1} \delta$, namely for 
any $\phi \in \mcal{K}_{X}^{q}$ we have
\[ \phi \cdot a = \phi \cdot \partial^1(1) = 
(-1)^{q + 1} \delta(\phi) . \]
\end{exa}

\begin{exa} \label{exa1.3}
If $X$ is a nonsingular curve then the theorem takes on a very 
simple form. Here 
$\mcal{K}^{-1}_{X} = \Omega^{1}_{K / \k}$,
where $K$ is the function field of $X$; and
$\mcal{K}^{0}_{X} = \opn{Coker}(\Omega^{1}_{X / \k} \to 
\Omega^{1}_{K / \k})$.
And there is actually an isomorphism of complexes
\[ \ul{\mbb{A}}_{\mrm{red}}^{\bdot}(\mcal{O}_{X}) \cong 
\mcal{H}om_{\mcal{O}_{X}}
(\mcal{K}^{\bdot}_{X}, \mcal{K}^{\bdot}_{X}) \]
(cf.\ Example \ref{exa1.2}).
\end{exa}

If $X$ is integral of dimension $n$, 
let $\bsym{\omega}_{X}$ be the coherent sheaf
$\mrm{H}^{-n} \mcal{K}_{X}^{\bdot}$. This is a subsheaf of the constant 
sheaf 
$\mcal{K}^{-n}_{X} = \Omega^{n}_{K / \k}$,
where $K$ is the function field of $X$.
In \cite{Li} and \cite{KW}
$\bsym{\omega}_{X}$ is called the sheaf of regular differential forms
(cf.\ \cite{Ye1} Theorem 4.4.16).
Since $\mcal{K}^{\bdot}_{X}$ is a complex of injectives and
$\bsym{\omega}_{X} \to \ul{\mbb{A}}_{\mrm{red}}^{\bdot}
(\bsym{\omega}_{X})$
is a quasi-isomorphism, there is a map of complexes
$\ul{\mbb{A}}_{\mrm{red}}^{\bdot}(\bsym{\omega}_{X}) 
\to \mcal{K}^{\bdot}_{X}[-n]$
inducing the identity in $\mrm{H}^0$.
Lipman asked for an explicit formula for such a homomorphism 
of complexes. Producing such a formula was the main result 
of \cite{HY1}. The following corollary gives essentially 
the same formula but in terms of DG modules.

\begin{cor} \label{cor1.1}
Suppose $X$ is integral of dimension $n$. Then the map
\[ \ul{\mbb{A}}_{\mrm{red}}^{\bdot}(\bsym{\omega}_{X}) \cong
\bsym{\omega}_{X} \otimes_{\mcal{O}_{X}} 
\ul{\mbb{A}}_{\mrm{red}}^{\bdot}(\mcal{O}_{X}) \to
\mcal{K}_{X}^{\bdot}[-n] , \]
sending
$\alpha \otimes a \mapsto  \alpha \cdot a$,
is a homomorphism of complexes,
which induces the identity on
$\bsym{\omega}_{X} = 
\mrm{H}^{0} \ul{\mbb{A}}_{\mrm{red}}^{\bdot}(\bsym{\omega}_{X}) =
\mrm{H}^{0} \mcal{K}_{X}^{\bdot}[-n]$.
\end{cor}

\begin{proof}
In 
$\bsym{\omega}_{X} \otimes_{\mcal{O}_{X}} 
\ul{\mbb{A}}_{\mrm{red}}^{\bdot}(\mcal{O}_{X})$
one has 
$\partial(\alpha \otimes a) = \alpha \otimes \partial(a)$.
On the other hand in $\mcal{K}_{X}^{\bdot}[-n]$ the section
$\alpha$ has degree $0$, and $\delta(\alpha) = 0$, so 
\[ \delta (\alpha \cdot a) = \delta(\alpha) \cdot a + 
\alpha \cdot \partial(a) = \alpha \cdot \partial(a) . \]

\end{proof}

Given a morphism $f : X \to Y$ there is a natural DGA homomorphism
$f^{*} : \ul{\mbb{A}}_{\mrm{red}}^{\bdot}(\mcal{O}_{Y}) \to
f_{*} \ul{\mbb{A}}_{\mrm{red}}^{\bdot}(\mcal{O}_{X})$.
There is also a map of graded $\mcal{O}_{Y}$-modules
$\opn{Tr}_{f} : f_{*} \mcal{K}^{\bdot}_{X} \to 
\mcal{K}^{\bdot}_{Y}$,
and this is a map of complexes if $f$ is proper (see \cite{Ye3} 
Definition 2.11 and Theorem 3.4).

\begin{thm} \label{thm1.2}
If $f : X \to Y$ is proper then
$\opn{Tr}_{f} : f_{*} \mcal{K}^{\bdot}_{X} \to \mcal{K}^{\bdot}_{Y}$
is a homomorphisms of DG 
$\ul{\mbb{A}}_{\mrm{red}}^{\bdot}(\mcal{O}_{Y})$-modules.
\end{thm}

\begin{proof}
We have to show that 
$\opn{Tr}_{f}(\phi \cdot a) = \opn{Tr}_{f}(\phi) \cdot a$
for local sections $\phi \in \mcal{K}^{\bdot}_{X}$ and 
$a \in \ul{\mbb{A}}_{\mrm{red}}^{\bdot}(\mcal{O}_{Y})$.
For this we may as well assume 
$\phi = \phi_{x_0} \in \mcal{K}_{X}(x_{0})$
and $a = a_{\eta} \in \mcal{O}_{Y, \eta}$ 
for a saturated chain 
$\eta = (y_{0}, \ldots, y_{q})$ in $Y$
and a point $x_{0}$ which is closed in the fiber
$f^{-1}(y_{0})$. According to 
\cite{Ye3} Proposition 2.1 we have an isomorphism of BCAs
\[ \prod_{\xi} \mcal{O}_{X, \xi} \cong
\mcal{O}_{X, (x_{0})} \otimes^{(\wedge)}_{\mcal{O}_{Y, (y_{0})}}
\mcal{O}_{Y, \eta}  \]
where $\xi = (x_{0}, \ldots, x_{q})$ runs over the 
finitely many chains in $X$ satisfying $f(x_{i}) = y_{i}$.
Therefore by \cite{Ye2} Theorem 7.4 the left square in the diagram
\[ \begin{CD}
\mcal{K}(\mcal{O}_{X, (x_{0})}) @>\mrm{q}>>
\bigoplus_{\xi} \mcal{K}(\mcal{O}_{X, \xi}) @>\opn{Tr}>> 
\bigoplus_{x_{q}} \mcal{K}(\mcal{O}_{X, (x_{q})}) \\
@V\opn{Tr}VV @V\opn{Tr}VV @V\opn{Tr}VV \\
\mcal{K}(\mcal{O}_{Y, (y_{0})}) @>\mrm{q}>>
\mcal{K}(\mcal{O}_{Y, \eta})@> \opn{Tr}>> 
\mcal{K}(\mcal{O}_{Y, (y_{q})}) \\
\end{CD} \] 
is commutative. The functoriality of $\opn{Tr}$ with respect to 
morphisms of BCAs implies the commutativity of the right square.
\end{proof}

\begin{que} \label{que1.1}
It is plausible to assume that an action as in Theorem 
\ref{thm1.1} exists even without the explicit construction of the 
residue complex. Suppose $X$ is a finite dimensional
noetherian scheme endowed with a dimension function 
$d : X \to \mbb{Z}$ (e.g.\ if $X$ has a dualizing complex). Let 
$\mcal{L}$ be a quasi-coherent Cousin complex on $X$, so that 
$\mcal{L}^{-q} \cong \boplus_{d(x) = q} \mrm{H}^{-q}_x \mcal{L}$.
Is $\mcal{L}$ a right DG 
$\ul{\mbb{A}}_{\mrm{red}}^{\bdot}(\mcal{O}_{X})$-module?
\end{que}


\section{de Rham Complexes, Residues and Adeles}

Let us continue with the setup of Section 1.
Consider the algebraic de Rham complex
$(\Omega^{\bdot}_{X/\k}, \mrm{d})$ of $X$. 
For any integers $p,q$ define
\[ \mcal{A}^{p, q}_{X} := \ul{\mbb{A}}_{\mrm{red}}^{q}(
\Omega^{p}_{X / \k}) . \]
Since $\mcal{M} \mapsto \ul{\mbb{A}}_{\mrm{red}}^{q}(\mcal{M})$
is functorial with respect to differential operators
$\mcal{M} \to \mcal{N}$ (cf.\ \cite{HY1}), the sheaves 
$\mcal{A}^{p, q}_{X}$
make up a double complex, with commuting operators
$\mrm{d}$ and $\partial$. Let
\[ \mcal{A}_{X}^{n} := \bigoplus_{p+q=n} \mcal{A}_{X}^{p,q} \]
with coboundary operators
\[ \begin{split}
\mrm{D}' & := \mrm{d} : \mcal{A}_{X}^{p,q} \to \mcal{A}_{X}^{p+1,q} , \\ 
\mrm{D}'' & := (-1)^{p} \partial : 
\mcal{A}_{X}^{p,q} \to \mcal{A}_{X}^{p,q+1} 
\end{split} \]
and 
$\mrm{D} := \mrm{D}' + \mrm{D}''$.
The Alexander-Whitney product
$\mcal{A}_{X}^{p,q} \otimes \mcal{A}_{X}^{p',q'} \to
\mcal{A}_{X}^{p+p',q+q'}$ 
makes $(\mcal{A}_{X}^{\bdot}, \mrm{D})$ into a DGA. 
To be explicit, write
$\mcal{A}_{X}^{p,q} \cong 
\Omega^{p}_{X/\k} \otimes_{\mcal{O}_{X}}
\ul{\mbb{A}}_{\mrm{red}}^{q}(\mcal{O}_{X})$. 
Then taking local sections
$a  \in \ul{\mbb{A}}_{\mrm{red}}^{q}(\mcal{O}_{X})$,
$b  \in \ul{\mbb{A}}_{\mrm{red}}^{q'}(\mcal{O}_{X})$,
$\alpha \in \Omega^{p}_{X/\k}$ and
$\beta \in \Omega^{p'}_{X/\k}$, one has
\[ (\alpha \otimes a) \cdot (\beta \otimes b) =
(-1)^{qp'} \alpha \wedge \beta \otimes 
\partial^{-} (a) \cdot
\partial^{+} (b) \in
\mcal{A}_{X}^{p+p',q+q'}  \]
and
\[ \begin{aligned}
\mrm{D}(\alpha \otimes a) & = 
\mrm{D}(\alpha) \cdot a + (-1)^p \alpha \cdot \mrm{D}(a) \\
& = \mrm{d}(\alpha) \cdot a + (-1)^p \alpha \cdot \mrm{d}(a) 
+ (-1)^p \alpha \cdot \partial(a) . 
\end{aligned} \]
Each $\mcal{A}_{X}^{p,q}$ is a flasque sheaf, and the natural 
homomorphism of DGAs
$\Omega^{\bdot}_{X / \k} \to \mcal{A}_{X}^{\bdot}$ is a 
quasi-isomorphism.

In \cite{Ye3} it is proved that any differential operator
$D : \mcal{M} \to \mcal{N}$
between $\mcal{O}_{X}$-modules induces a dual operator
\[ \opn{Dual}(D) : \opn{Dual} \mcal{N} \to \opn{Dual} \mcal{M} , 
\]
where by definition
$\opn{Dual} \mcal{M}$ is the complex
$\mcal{H}om_{\mcal{O}_{X}} (\mcal{M}, \mcal{K}^{\bdot}_{X})$.
The operator \linebreak
$\opn{Dual}(D)$ is defined locally, 
in terms of Beilinson completion algebras. Let
\[ \mcal{F}_{X}^{p, q} := \mcal{H}om_{X}
(\Omega^{-p}_{X / \k}, \mcal{K}^{q}_{X}) . \]
We see that $\mcal{F}_{X}^{\bdot, \bdot}$ 
is a double complex, with two commuting operators 
$\delta$ and $\opn{Dual}(\mrm{d})$.
Set
\begin{equation} \label{eqn2.2}
\begin{aligned}
\mcal{F}_{X}^{n} & := \bigoplus_{p+q = n} \mcal{F}_{X}^{p, q} , \\
\mrm{D}' & := (-1)^{p + q + 1} \opn{Dual}(\mrm{d}) : 
\mcal{F}_{X}^{p, q} \to \mcal{F}_{X}^{p+1, q} , \\
\mrm{D}'' & := \delta : \mcal{F}_{X}^{p, q} \to \mcal{F}_{X}^{p, q+1}
\quad \text{and} \\
\mrm{D} & := \mrm{D}' + \mrm{D}'' : 
\mcal{F}_{X}^{n} \to \mcal{F}_{X}^{n + 1} .
\end{aligned}
\end{equation}

\begin{proof}[Proof of Theorem \tup{\ref{thm0.2}}]
The left action of $\Omega^{\bdot}_{X/\k}$ on itself makes 
$\mcal{F}^{\bdot}_{X}$ into a graded right 
$\Omega^{\bdot}_{X/\k}$-module. 
By Theorem \ref{thm1.1}, $\mcal{F}^{\bdot}_{X}$ is a graded right 
$\ul{\mbb{A}}^{\bdot}_{\mrm{red}}(\mcal{O}_{X})$-module. The formula 
for the product $\phi \cdot a$ which makes the signs correct is 
\[ (\phi \cdot a)(\beta) := (-1)^{pq'} \phi(\beta) \cdot a \]
for local sections 
$\phi \in \mcal{F}^{p,q}_{X}$,  
$a \in \ul{\mbb{A}}^{q'}_{\mrm{red}}(\mcal{O}_{X})$
and
$\beta \in \Omega^{-p}_{X/\k}$.
Since 
$\mcal{A}^{\bdot}_{X} \cong 
\Omega^{\bdot}_{X/\k} \otimes_{\mcal{O}_{X}}
\ul{\mbb{A}}^{\bdot}_{\mrm{red}}(\mcal{O}_{X})$
as graded algebras, we obtain a structure of graded right 
$\mcal{A}^{\bdot}_{X}$-module on $\mcal{F}^{\bdot}_{X}$.

It remains to check the coboundaries, which we break up into four 
steps, calculating 
$\mrm{D}'(\phi \cdot a)$, $\mrm{D}'(\phi \cdot \alpha)$,
$\mrm{D}''(\phi \cdot a)$ and $\mrm{D}''(\phi \cdot \alpha)$
separately, with $a, \phi$ be as above, and  
$\alpha \in \Omega^{p'}_{X/\k}$.

Since $\mrm{D}''(\alpha) = 0$ we get
\[ \begin{aligned}
\mrm{D}''(\phi \cdot \alpha)(\beta) 
& = \delta((\phi \cdot \alpha)(\beta))
= \delta(\phi(\alpha \wedge \beta)) \\
& = \mrm{D}''(\phi)(\alpha \wedge \beta)
= (\mrm{D}''(\phi) \cdot \alpha)(\beta) \\
& = (\mrm{D}''(\phi) \cdot \alpha + (-1)^{p + q} \phi \cdot 
\mrm{D}''(\alpha))(\beta) 
\end{aligned} \]
for every 
$\beta \in \Omega^{-p}_{X / \k}$.

Next, by Thm.\ \ref{thm1.1}
\[ \begin{aligned}
(\mrm{D}''(\phi \cdot a))(\beta) 
& = \delta((\phi \cdot a)(\beta)) 
= \delta((-1)^{p q'}\phi(\beta) \cdot a) \\
& = (-1)^{pq'} \mrm{D}''(\phi(\beta) \cdot a) \\
& = (-1)^{pq'} \left( \mrm{D}'' (\phi(\beta)) \cdot a + 
(-1)^{q} \phi(\beta) \cdot \mrm{D}''(a) \right) \\
& = \left( \mrm{D}''(\phi) \cdot a + (-1)^{p+q} \phi \cdot 
\mrm{D}''(a) \right)(\beta) . 
\end{aligned} \]
This takes care of $\mrm{D}''$. 

As for $\mrm{D}'$, by \cite{Ye3} Proposition 5.2 we have
\begin{equation} \label{eqn2.3}
\mrm{D}' (\phi \cdot \alpha) = \mrm{D}'(\phi) \cdot \alpha +
(-1)^{p+q} \phi \cdot \mrm{D}'(\alpha) . 
\end{equation}

Finally we will prove that
\[ \mrm{D}'(\phi \cdot a) = 
\mrm{D}'(\phi) \cdot a + (-1)^{p + q} \phi \cdot \mrm{D}'(a) \]
by induction on $q'$. As before we may assume
\[ \phi = \phi_x \in \mcal{H}om_{\mcal{O}_{X}}
(\Omega^{-p}_{X / \k}, \mcal{K}_X(x)) \subset
\mcal{F}_X^{p, q} \]
for some point $x = x_0$, and
\[ a = a_{\xi} \in \mcal{O}_{X, \xi} \subset 
\ul{\mbb{A}}^{q'}_{\mrm{red}}(\mcal{O}_{X}) \]
for some saturated chain
$\xi = (x_0, \ldots, x_{q'})$.
Choose an integer $n$ such that
$\mfrak{m}_x^n \cdot \phi_x = 0$. 
Then we also have
$\mfrak{m}_x^{n + 1} \cdot \mrm{D}'(\phi_x) = 0$. 

For $q' = 0$ we have $a_{(x)} \in \mcal{O}_{X, (x)}$. Choose some 
local section $b \in \mcal{O}_{X}$ near $x$ such that 
$b \equiv a_{(x)} \opn{mod} \mfrak{m}_x^{n + 1}$.
So by Lemma \ref{lem1.4} one has
$\phi \cdot a = \phi \cdot b$,
$\mrm{D}'(\phi) \cdot a = \mrm{D}'(\phi) \cdot b$
and 
$\phi \cdot \mrm{D}'(a) = \phi \cdot \mrm{D}'(b)$. 
But by equation (\ref{eqn2.3}), with 
$\alpha = b \in \Omega^0_{X / \k}$, we know that
\begin{equation} \label{eqn2.4}
\mrm{D}'(\phi \cdot b) = \mrm{D}'(\phi) \cdot b +
(-1)^{p+q} \phi \cdot \mrm{D}'(b) . 
\end{equation}

Now let us handle the case $q' = 1$, i.e.\ 
$\xi = (x_0, x_1)$. By Lemma \ref{lem1.4}
there exist adeles
$b \in \mcal{O}_{X, (x_0)}$
and
$c \in \mcal{O}_{X, (x_1)}$
such that
\[ a_{\xi} = b_{(x_0, x_1)} \cdot c_{(x_1)}
\in \mcal{O}_{X, \xi} /  \mfrak{m}_{x_0}^{n + 1} \mcal{O}_{X, \xi} 
, \]
where $b_{(x_0, x_1)}$ is the $(x_0, x_1)$ component of 
$\partial^1(b)$. Let 
$d := 1_{(x_0, x_1)} \in \mcal{O}_{X, (x_0, x_1)} \subset
\ul{\mbb{A}}^{1}_{\mrm{red}}(\mcal{O}_{X})$.
We then get
$\partial^1(b) = b \cdot d$, 
\[ a - b \cdot d \cdot c
\in \mfrak{m}_{x_0}^{n + 1} \mcal{O}_{X, \xi} , \]
$\phi \cdot a = \phi \cdot (b \cdot d \cdot c)$, 
$\mrm{D}'(\phi) \cdot a = \mrm{D}'(\phi) \cdot (b \cdot d \cdot c)$
and
$\phi \cdot \mrm{D}'(a) = \phi \cdot \mrm{D}'(b \cdot d \cdot c)$.
By the $q' = 0$ case we know that equation (\ref{eqn2.4}) holds, 
and also that
\[ \mrm{D}'((\phi \cdot b \cdot d) \cdot c) 
= \mrm{D}'(\phi  \cdot b \cdot d) \cdot c +
(-1)^{p + q + 1} (\phi \cdot b \cdot d) \cdot \mrm{D}'(c) . \]
It remains to verify that
\[ \mrm{D}'(\psi \cdot d) 
= \mrm{D}'(\psi) \cdot d +
(-1)^{p + q} \psi \cdot \mrm{D}'(d) , \]
where we define $\psi := \phi \cdot b$. But for every 
$\beta \in \Omega^{-p}_{X / \k}$ we have
\[ (\psi \cdot d)(\beta) = (-1)^p \psi(\beta) \cdot d = 
(-1)^p \delta_{(x_0, x_1)}(\psi(\beta)) = 
(-1)^{p + q + 1} \mrm{D}''(\psi)(\beta) . \]
Therefore 
$\psi \cdot d = (-1)^{p + q + 1} \mrm{D}''(\psi)$.
Since $\mrm{D}' \circ \mrm{D}'' = - \mrm{D}'' \circ \mrm{D}'$
and
$\mrm{D}'(d) = 0$
we arrive at
\[ \begin{aligned}
\mrm{D}'(\psi \cdot d) 
& = (-1)^{p + q + 1} (\mrm{D}' \circ \mrm{D}'')(\psi)
= (-1)^{p + q} (\mrm{D}'' \circ \mrm{D}')(\psi) \\
& = (-1)^{p + q} \mrm{D}''(\mrm{D}'(\psi))
= (-1)^{(p + q) + (p + q + 2)} \mrm{D}'(\psi) \cdot d \\
& = \mrm{D}'(\psi) \cdot d + 
(-1)^{p + q} \psi \cdot \mrm{D}'(d) . 
\end{aligned} \]

To conclude we consider $q' \geq 2$. Since we are working locally, 
we can assume by the approximation trick (Lemma \ref{lem1.4}) that 
$a = b \cdot c$ with 
$b \in \ul{\mbb{A}}^{1}_{\mrm{red}}(\mcal{O}_{X})$
and
$c \in \ul{\mbb{A}}^{q' - 1}_{\mrm{red}}(\mcal{O}_{X})$. 
The induction hypothesis applies to $b$ and $c$, and we have
\[ \begin{aligned}
\mrm{D}'(\phi \cdot a) 
& = \mrm{D}'(\phi \cdot b \cdot c) \\
& = \mrm{D}'(\phi \cdot b) \cdot c + (-1)^{p + q + 1} 
\phi \cdot b \cdot \mrm{D}'(c) \\
& = \mrm{D}'(\phi) \cdot b \cdot c + 
(-1)^{p + q} \phi \cdot \mrm{D}'(b) \cdot c +
(-1)^{p + q + 1} 
\phi \cdot b \cdot \mrm{D}'(c) \\
& = \mrm{D}'(\phi) \cdot b \cdot c + 
(-1)^{p + q} \phi \cdot \mrm{D}'(b \cdot c) \\
& = \mrm{D}'(\phi) \cdot a + 
(-1)^{p + q} \phi \cdot \mrm{D}'(a) .
\end{aligned} \]
\end{proof}

\begin{proof}[Proof of Corollary \tup{\ref{cor0.1}}]
By \cite{Ye3} Proposition 5.4 the trace homomorphisms
$\opn{Tr}: \mcal{K}_{X_n}^{\bdot} \to 
\mcal{K}_{X_{n + 1}}^{\bdot}$
give rise to DG $\Omega^{\bdot}_{X_{n + 1} / \k}$-linear 
homomorphisms 
$\opn{Tr}: \mcal{F}_{X_n}^{\bdot} \to 
\mcal{F}_{X_{n + 1}}^{\bdot}$. 
Therefore $\mcal{F}^{\bdot}_{\mfrak{X}}$ is a right DG 
$\mcal{A}^{\bdot}_{\mfrak{X}}$-module. 

Let $\widehat{\Omega}^{\bdot}_{\mfrak{X} / \k}$ be the formal 
completion of $\Omega^{\bdot}_{Y / \k}$ along $X$. According to 
\cite{HY1} the DGA quasi-isomorphism
$\widehat{\Omega}^{\bdot}_{\mfrak{X} / \k} \to
\mcal{A}^{\bdot}_{\mfrak{X}}$
is a flasque resolution, so
\[ \mrm{H}^i(X, \widehat{\Omega}^{\bdot}_{\mfrak{X} / \k})
= \mrm{H}^i \Gamma(X, \mcal{A}^{\bdot}_{\mfrak{X}}) . \] 
By \cite{Ye3} Proposition 6.8 we get 
$\mcal{F}^{\bdot}_{\mfrak{X}} \cong
\mrm{R} \Gamma_X \Omega^{\bdot}_{Y / \k}[2m]$,
where $m := \opn{dim} Y$, and hence
\[ \mrm{H}^{2m - i}_X(Y, \Omega^{\bdot}_{Y / \k}) = 
\mrm{H}^{-i} \Gamma(X, \mcal{F}^{\bdot}_{\mfrak{X}}) , \]
compatible with the cap product. But when $\opn{char} \k = 0$ 
one has by definition
$\mrm{H}^i_{\mrm{DR}}(X) :=  
\mrm{H}^i(X, \widehat{\Omega}^{\bdot}_{\mfrak{X} / \k})$
and
$\mrm{H}_i^{\mrm{DR}}(X) := 
\mrm{H}^{2m - i}_X(Y, \Omega^{\bdot}_{Y / \k})$,
see \cite{Ha}.
\end{proof}

\begin{rem} \label{rem2.2}
The role of the residue complex as a canonical flasque resolution of 
$\Omega^{\bdot}_{X/\k}$ (for $X$ smooth) appears already in \cite{Ha}
and in \cite{EZ}. In order to define the operator dual to $\mrm{d}$
these authors use the fact that 
$\mcal{F}^{\bdot}_{X} \cong
\Omega^{\bdot}_{X/\k} \otimes_{\mcal{O}_{X}} \mcal{K}^{\bdot}_{X}$
is the Cousin resolution of $\Omega^{\bdot}_{X/\k}$ in
the category of abelian sheaves on $X$. This should be compared to 
our approach (in \cite{Ye3} and here) where the 
dual operator $\opn{Dual}(\mrm{d})$
is defined in terms of BCAs, i.e.\ by algebraic-analytic methods.
\end{rem}

\begin{rem} \label{rem2.1}
Suppose $X$ is smooth irreducible of dimension $n$. One can regard 
the sheaf 
$\mcal{A}_{X}^{p, q}$ as an analog of the Dolbeault sheaf of smooth
$(p, q)$ forms on a complex manifold (cf.\ \cite{GH} Ch.\ 0). The 
adelic resolution
\[ 0 \to \Omega^{p}_{X / \k} \to \mcal{A}_{X}^{p, 0} 
\xrightarrow{\mrm{D}''} \mcal{A}_{X}^{p, 1} \to \cdots \to
\mcal{A}_{X}^{p, n} \to 0 \]
corresponds to the $\bar{\partial}$ resolution on the manifold. 
Any section 
$\phi \in \Gamma(X, \mcal{F}_{X}^{-p, -q})$ 
determines a functional on $\Gamma(X, \mcal{A}_{X}^{p, q})$, namely
$\alpha \mapsto \opn{Tr}_{X / \k} (\phi \cdot \alpha) \in \k$.
Here 
\[ \opn{Tr}_{X / \k} : \Gamma(X, \mcal{F}_{X}^{0, 0}) 
= \Gamma(X, \mcal{K}_{X}^{0}) \to \k \]
is the ``sum of residues'' (cf.\ \cite{Ye3} Def.\ 1.16). 
In this way we can think of $\mcal{F}_{X}^{-p, -q}$ as an analog of 
the sheaf of $(p, q)$ currents on a manifold (cf.\ \cite{GH} Ch.\ 3).
\end{rem}

\begin{rem} \label{rem2.3}
Suppose $X$ is smooth irreducible and $\opn{char} \k = 0$.
Let $Z \subset X$ be an irreducible closed subset of codimension $d$,
and let
\[ 0 \to \mcal{E}_{l} \to \cdots \to \mcal{E}_{1} \to 
\mcal{E}_{0} \to \mcal{O}_{Z} \to 0 \]
be a finite locally free resolution. Using the adelic Chern-Weil
theory of \cite{HY2} it is possible to construct a Chern character
form 
$\opn{ch}(\mcal{E}; \nabla) \in \mcal{A}^{\bdot}_{X}$, 
depending on adelic connections $\nabla_{i}$ on the $\mcal{E}_{i}$, 
whose component
$\opn{ch}(\mcal{E}; \nabla)_{2d} \in \mcal{A}^{2d}_{X}$ 
satisfies
\[ \mrm{C}_{X} \cdot \opn{ch}(\mcal{E}; \nabla)_{2d} = \mrm{C}_{Z}
\in \mcal{F}^{\bdot}_{X} . \]
Here $\mrm{C}_{X} \in \mcal{F}^{-2n}_{X}$ and
$\mrm{C}_{Z} \in \mcal{F}^{-2(n - d)}_{X}$
are the fundamental classes.
This generalizes \cite{HY2} Theorem 6.5. 
\end{rem}


\end{document}